\documentclass{elsart3-1}

\usepackage{amsmath,amssymb,graphicx,epsfig,color}
\usepackage[english,francais]{babel}

\newtheorem{e-proposition}[theorem]{Proposition}

\newtheorem{e-definition}[theorem]{Definition\rm}

\newtheorem{theoreme}{Th\'eor\`eme}[section]

\newtheorem{proposition}[theoreme]{Proposition}

\renewcommand{\theequation}{\arabic{equation}}
\def\zz{\mathbb{Z}}
\def\rr{\mathbb{R}}

\def\og{\leavevmode\raise.3ex\hbox{$\scriptscriptstyle\langle\!\langle$~}}
\def\fg{\leavevmode\raise.3ex\hbox{~$\!\scriptscriptstyle\,\rangle\!\rangle$}}

\begin{document}

\begin{frontmatter}

\selectlanguage{english}
\title{Localized solutions and filtering mechanisms for the discontinuous Galerkin semi-discretizations of the $1-d$ wave equation}

\vspace{-2.6cm}

\selectlanguage{francais}
\title{Solutions localis\'{e}es et m\'{e}canismes de filtrage pour les approximations Galerkin discontinu de l'\'{e}quation des ondes}

\selectlanguage{english}
\author[BCAM]{Aurora Marica},
\ead{marica@bcamath.org}
\author[Iker,BCAM]{Enrique Zuazua}
\ead{zuazua@bcamath.org}

\address[Iker]{Ikerbasque, Basque Foundation for Science, Alameda Urquijo 36-5, Plaza Bizkaia, 48011, Bilbao, Basque Country, Spain}

\address[BCAM]{BCAM - Basque Center for Applied Mathematics, Bizkaia Technology Park 500, 48160, Derio, Basque Country, Spain}

\begin{abstract}
We perform a complete Fourier analysis of the semi-discrete  $1-d$ wave equation obtained through a $P_1$
discontinuous Galerkin (DG) approximation of the continuous wave equation on an uniform grid. The resulting system exhibits the interaction of two types of components:
a physical one and a spurious one, related to the
possible discontinuities that the numerical solution allows. Each dispersion relation contains critical points where the corresponding group velocity vanishes. Following previous constructions, we rigorously build wave packets with arbitrarily small velocity of propagation concentrated either on the physical or on the spurious component. We also develop filtering mechanisms aimed at recovering the uniform velocity of propagation of the continuous solutions. Finally, some applications to numerical approximation issues of control problems are also presented.

\vskip 0.5\baselineskip

\selectlanguage{francais}
\noindent{\bf R\'esum\'e
} \vskip 0.5\baselineskip \noindent
On d\'{e}veloppe une analyse de Fourier compl\`{e}te de l'\'{e}quation des ondes unidimensionnelle semi-discr\'{e}tis\'{e}e en espace obtenue dans l'approximation num\'{e}rique de l'\'{e}quation des ondes par une m\'{e}thode de Galerkin discontinue (GD) $P_1$ dans un maillage uniforme. On met en \'{e}vidence la coexistence de deux composantes dans le syst\`{e}me num\'{e}rique: une physique, et une parasite li\'{e}e aux discontinuit\'{e}s que la solution num\'{e}rique permet. Chaque relation de dispersion contient des points critiques o\`{u} la vitesse de groupe correspondante s'annule. En suivant les constructions faites ant\'{e}rieurement pour le sch\'{e}ma en diff\'{e}rences finies, on construit d'une mani\`{e}re rigoureuse des paquets d'ondes qui se propagent \`{a} une vitesse arbitrairement petite, concentr\'{e}s soit sur une composante ou sur l'autre. On d\'{e}veloppe aussi des m\'{e}canismes de filtrage permettant de r\'{e}cup\'{e}rer les propri\'{e}t\'{e}s de propagation des solutions de l'\'{e}quation continue. Enfin, on pr\'{e}sente une application \`{a} l'approximation num\'{e}rique des probl\`{e}mes de contr\^{o}le.
\end{abstract}
\end{frontmatter}

\selectlanguage{francais}
\section*{Version fran\c{c}aise abr\'eg\'ee}

Dans cet article, on consid\`{e}re le probl\`{e}me de Cauchy associ\'{e} \`{a} l'\'{e}quation des ondes unidimensionnelle semi-discretis\'{e}e en espace obtenue dans l'approximation num\'{e}rique de l'\'{e}quation des ondes continue par une m\'{e}thode de Galerkin discontinue (GD) appel\'{e}e Symmetric Interior Penalty Discontinuous Galerkin (SIPG) (cf. \cite{BreUnif}) dans un maillage uniforme de pas $h$ en utilisant des polyn\^{o}mes du premier ordre. Ce sch\'{e}ma produit deux relations de dispersion: une physique, associ\'{e}e \`{a} la partie continue de la solution num\'{e}rique, et une parasite, associ\'{e}e \`{a} ses sauts. On d\'{e}veloppe d'abord l'analyse de Fourier de cette approximation qui met en \'{e}vidence le fait que la vitesse de groupe correspondante \`{a} chaque relation de dispersion s'annule pour certains nombres d'onde. Cela nous permet d'adapter la construction de paquets d'ondes ant\'{e}rieurement r\'{e}alis\'{e}e pour le sch\'{e}ma en diff\'{e}rences finies en \cite{ErvZuaBook}, \cite{MarZuaCRASFD} et de construire des solutions num\'{e}riques concentr\'{e}es sur l'un des deux modes possibles. Cela fait que la propri\'{e}t\'{e} d'observabilit\'{e} bien connue pour l'\'{e}quation des ondes continue, \`{a} savoir que pour un temps suffisamment grand l'\'{e}nergie totale des solutions peut \^{e}tre estim\'{e}e en fonction de l'\'{e}nergie localis\'{e}e dans le compl\'{e}mentaire d'un ensemble compact (cf. \cite{ZuaUnbDom}), n'est pas v\'{e}rifi\'{e}e de mani\`{e}re uniforme par rapport \`{a} $h$ dans ce cas semi-discret. On exhibe ainsi un autre exemple de ph\'{e}nom\`{e}ne pathologique concernant les propri\'{e}t\'{e}s de propagation et dispersion des approximations num\'{e}riques classiques de l'\'{e}quation d'ondes (\cite{ZuaPOC}, \cite{ErvZuaBook}, \cite{MarZuaCRASFD}) et de Schr\"{o}dinger (cf. \cite{LivEZDispSchr}).

Notre second objectif est de d\'{e}velopper des m\'{e}canismes de filtrage pour construire des classes de donn\'{e}es initiales dans lesquelles la propri\'{e}t\'{e} d'observabilit\'{e} aie lieu uniform\'{e}ment par rapport \`{a} $h$. Pour ceci, on consid\`{e}re des donn\'{e}es initiales avec des sauts nuls, assurant que  l'\'{e}nergie totale des solutions soit domin\'{e}e par l'\'{e}nergie des projections sur le mode physique, pour ensuite filtrer les hautes fr\'{e}quences par un algorithme bigrille (\cite{GloLioHeBook}, \cite{GloStokesWaves}, \cite{LivEZ}).

Nos r\'{e}sultats compl\`{e}tent la litt\'{e}rature existante sur les m\'{e}thodes de Galerkin discontinus et en particulier \cite{BuffaEig}, o\`{u} on d\'{e}montre que la plupart de ces approximations sont spectralement correctes, et \cite{AinsSisp}, o\`{u} les propri\'{e}t\'{e}s dispersives et dissipatives de la version $hp$ des semi-discr\'{e}tisations par GD de l'\'{e}quation des ondes sont analys\'{e}es.

\selectlanguage{english}

\smallskip\noindent\textbf{1. Fourier analysis of the $P_1$ discontinuous Galerkin approximations of the $1-d$ wave equation.} Based on discontinuous finite-element spaces, the discontinuous Galerkin (DG) methods can handle elements of various types and shapes, irregular non-matching grids and even varying polynomial order. A particular class of DG methods for elliptic and parabolic problems are the so-called interior penalty (IP) ones, where the continuity is weakly enforced across the element interfaces, by adding suitable bilinear forms, called numerical fluxes, to the classical variational formulations (see \cite{BreUnif}). This note provides a further contribution to the analysis of DG methods. Among the existing and related works, we refer to \cite{BuffaEig}, where the eigenproblem associated to the Laplace operator discretized by means of DG methods is analyzed, showing that several DG methods provide spectrally correct approximations of the Laplace operator, and to \cite{AinsSisp}, where the dispersive and dissipative properties of $hp$ versions of various DG methods are studied.

In this paper, we deal with the simplest setting of the Symmetric Interior Penalty Discontinuous Galerkin (SIPG) (cf. \cite{BreUnif}) space semi-discretization of the $1-d$ wave equation on an uniform grid $x_j$, $j\in\zz$, with first order polynomials.
\begin{figure}[!]
  \begin{center}\includegraphics[width=6cm, height=4cm]{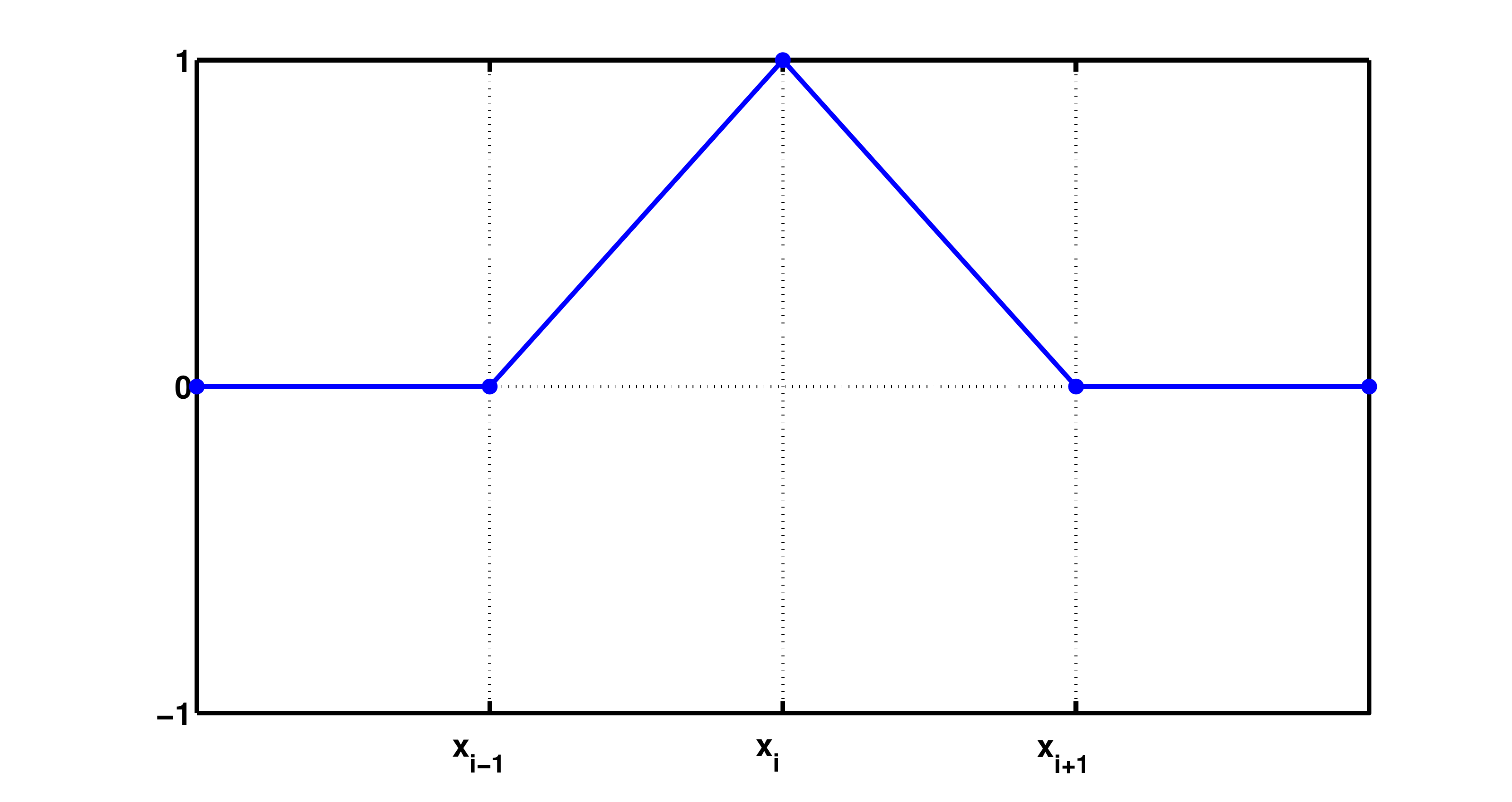}\includegraphics[width=6cm, height=4cm]{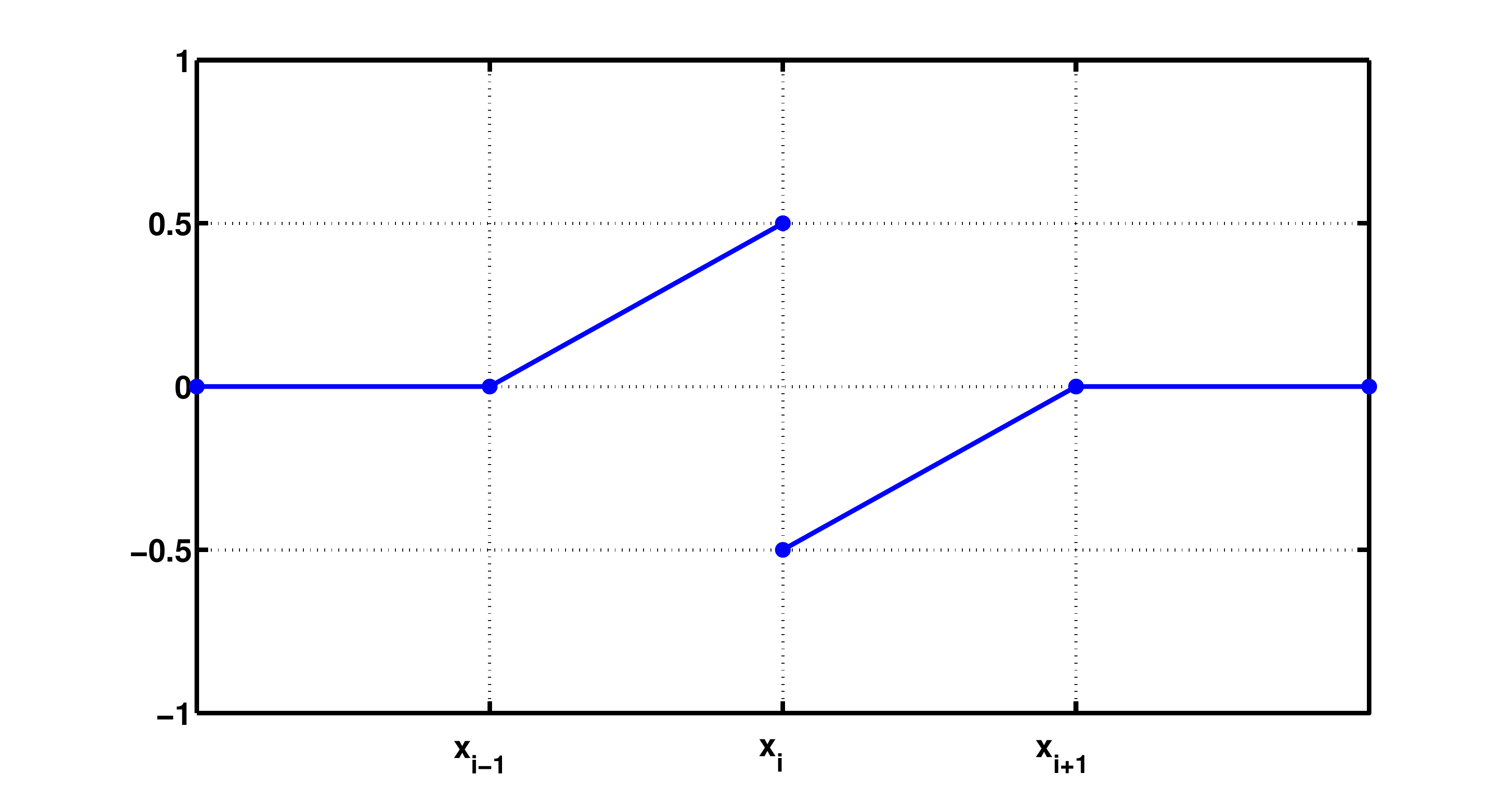}\\
  \caption{Basis functions for the $P_1$-DG methods: $\phi_i^{A}$ (left) and $\phi_i^{J}$ (right).}\end{center}\label{figbasis}
\end{figure}
The numerical solutions given by the DG methods being discontinuous, their numerical representation in $1-d$ consists in two values (possibly different ones) on every nodal point, representing the values to the left and to the right. An alternative, an often more convenient, way of representing the numerical solutions in the DG methods is in terms of the averages (denoted $\{\cdot\}$) and the jumps (denoted $[\cdot]$) along the interface, defined as $\{f\}(x)=(f(x+)+f(x-))/2$ and $[f](x)=f(x-)-f(x+)$. The finite element space is given by $V_h=U_h^{A}\oplus U_h^{J}$, with
$U_h^{A}=\mbox{span}\{\phi_i^{A},i\in\mathbf{Z}\}\cap H^1(\rr)$ and
$U_h^{J}=\mbox{span}\{\phi_i^{J},i\in\mathbf{Z}\}$,
where $$\phi_i^{A}(x)=\Big[1-\frac{|x-x_i|}{h}\Big]^+, \quad \phi_i^{J}(x)=\frac{1}{2}\mbox{sign}(x_i-x)\Big[1-\frac{|x-x_i|}{h}\Big]^+.$$
Observe that $\phi_i^{A}$ are the typical basis functions
of the $P_1$-classical finite element method, whereas
$\phi_i^{J}$ are
designed to represent the jumps at the nodal points. The functions $(\phi_i^{A},\phi_i^{J})_{i\in\zz}$ constitute a basis for $V_h$. In this way, any $f\in V_h$ can be uniquely represented as $f(x)=\sum_{i\in\mathbf{Z}}f_i^{A}\phi_i^{A}(x)+\sum_{i\in\mathbf{Z}}f_i^{J}\phi_i^{J}(x)
=f^{A}(x)+f^{J}(x)$, where $f^{A}$ and
$f^{J}$ are the continuous and the jump components of $f$,
respectively. In this way, the
piecewise linear discontinuous functions under consideration are perturbations of the classical piecewise linear and continuous ones
by jumps added at each nodal point.

Consider now the following Cauchy problem associated to the $1-d$ continuous wave equation:
\begin{equation}\partial_t^2u-\partial_x^2u=0,\ x\in\rr,\ t>0,\quad
u(x,0)=u^0(x),\ \partial_tu(x,0)=u^1(x),\ x\in\rr,\label{ContWaveEqn}\end{equation}
which is well posed for $(u^0,u^1)\in
\dot{H}^{1}(\rr)\times L^2(\rr)$, where $\dot{H}^1(\rr)$ denotes the completion of
$C_c^{\infty}(\rr)$ with respect to the semi-norm $\|\cdot\|_{\dot{H}^1(\rr)}=\|\partial_x\cdot\|_{L^2(\rr)}$.

For each value $s>1$ of the so-called penalty parameter, consider the symmetric bilinear form on $V_h\times V_h$
$$a_h^s(u,v)=\sum\limits_{j\in\mathbf{Z}}\Big[\int_{x_j}^{x_{j+1}}u_x(x)v_x(x)\,dx-
\big([u](x_j)\{v_x\}(x_j)+[v](x_j)\{u_x\}(x_j)\big)+\frac{s}{h}[u](x_j)[v](x_j)\Big]$$
and the following semi-discretization of the wave equation (\ref{ContWaveEqn}):
\begin{equation}u_h^s(x,t)\in V_h,\
\partial_t^2\int\limits_{\rr}u_h^s(x,t)v(x)\,dx+a_h^s(u_h^s(\cdot,
t),v)=0,\forall v\in V_h.\label{VarForm}\end{equation}

The above variational problem (\ref{VarForm}) is completed by numerical initial data of the form $u_h^s(x,0)=u_h^0(x)\in V_h$ and $\partial_tu_h^s(x,0)=u_h^1(x)\in V_h$. The unknown $u_h^s(x,t)$, being an element of $V_h$ for each
$t>0$, can be decomposed as $u_h^s(x,t)=\sum_{k\in \mathbf{Z}}A_k(t)\phi_k^{A}(x)+\sum_{k\in\mathbf{Z}}J_k(t)\phi_k^{J}(x)$. Denoting by $\overrightarrow{U}^h(t)=(A_k(t),J_k(t))_{k\in\zz}'$, the problem (\ref{VarForm}) can be written in the following matrix form ($'$ means transposition of a matrix):
\begin{equation}M_h\partial_t^2\overrightarrow{U}^h(t)+R_h^s\overrightarrow{U}^h(t)=0,\quad \overrightarrow{U}^h(0)
=\overrightarrow{U}^{h,0},\quad \partial_t\overrightarrow{U}^h(0)=\overrightarrow{U}^{h,1}\label{DiscWaveEqn},\end{equation}
where $M_h$ and $R_h^s$ are the block tri-diagonal symmetric mass and the stiffness matrices obtained by an infinite repetition of the stencils $m_h$, $r_h^s$, where $m_h$, $r_h^s$ are the matrices
\begin{equation}m_h=\left(
                      \begin{array}{cc|cc|cc}
                        \frac{h}{6} & -\frac{h}{12} & \frac{2h}{3} & 0 & \frac{h}{6} & \frac{h}{12}\smallskip\\
                        \frac{h}{12} & -\frac{h}{24} & 0 & \frac{h}{6} & -\frac{h}{12} & -\frac{h}{24} \\
                      \end{array}
                    \right),\ r_h^s=\left(
                      \begin{array}{cc|cc|cc}
                        -\frac{1}{h} & 0 & \frac{2}{h} & 0 & -\frac{1}{h} & 0\smallskip\\
                        0 & -\frac{1}{4h} & 0 & \frac{2s-1}{2h} & 0 & -\frac{1}{4h} \\
                      \end{array}
                    \right).
\nonumber\end{equation}

Set $\Pi_h=[-\pi/h,\pi/h]$. For $\xi\in\Pi_h$, let us denote by $\widehat{A}^h(\xi,t)$,
$\widehat{J}^h(\xi,t)$ the semi-discrete Fourier transforms (SDFT) of the sequences of averages, $\overrightarrow{A}^h(t)$, and of jumps, $\overrightarrow{J}^h(t)$ (cf. \cite{LivEZ}).  Similarly, by $\widehat{A}^{h,i}(\xi)$,
$\widehat{J}^{h,i}(\xi)$, $i=0,1$, we denote the SDFTs of the initial data $\overrightarrow{A}^{h,i}$, $\overrightarrow{J}^{h,i}$. Set $\widehat{U}^h(\xi,t):=(\widehat{A}^h(\xi,t),\widehat{J}^h(\xi,t))'$. The Fourier symbols of the mass and stiffness matrices are
\begin{equation}M_h(\xi)=\left(
                           \begin{array}{cc}
                             \frac{2+\cos(\xi h)}{3} & \frac{i\sin(\xi h)}{6} \\
                             -\frac{i\sin(\xi h)}{6} & \frac{2-\cos(\xi h)}{12} \\
                           \end{array}
                         \right), \quad R_h^s(\xi)=\left(
                                               \begin{array}{cc}
                                                 \frac{4}{h^2}\sin^2\big(\frac{\xi h}{2}\big) & 0 \\
                                                 0 & \frac{s-\cos^2\big(\frac{\xi h}{2}\big)}{h^2} \\
                                               \end{array}
                                             \right).
\nonumber\end{equation}
Let $S_h^s(\xi):=(M_h(\xi))^{-1}R_h^s(\xi)$, $(M_h(\xi))^{-1}$ being the inverse of the matrix $M_h(\xi)$. The system (\ref{DiscWaveEqn}) can be transformed into the following Cauchy problem associated to a system of two linear ODE's whose unknown is the vector function $\widehat{U}^h(\xi,t)$, depending on the frequency parameter $\xi$:
\begin{equation}\widehat{U}^h_{tt}(\xi,t)+S_h^s(\xi)\widehat{U}^h(\xi,t)=0,\ \xi\in\Pi_h,\ t>0,\quad\widehat{U}^h(\xi,0)=\widehat{U}^{h,0}(\xi),\quad \widehat{U}^h_t(\xi,0)=\widehat{U}^{h,1}(\xi),\ \xi\in\Pi_h.\label{DiscWaveFour}\end{equation}
Denote by $\Lambda_{ph,h}^s(\xi)$, $\Lambda_{sp,h}^s(\xi)$ the two eigenvalues of $S_h^s(\xi)$, by $\lambda_{ph,h}^s(\xi)$, $\lambda_{sp,h}^s(\xi)$ their square roots and by $P_{ph,h}^s(\xi)$, $P_{sp,h}^s(\xi)$ the two corresponding eigenvectors. The notation $ph$ and $sp$ stands for the "physical" and the "spurious" components, respectively. In this way, $\lambda_{ph,h}^s(\xi)$, $\lambda_{sp,h}^s(\xi)$ are the physical and the spurious dispersion relations. On the other hand, let $\omega_h(\xi)$, $\lambda_h(\xi)$ be the dispersion relations corresponding to the finite difference (FD) and $P_1$ classical finite element semi-discretizations of the $1-d$ wave equation. The physical dispersion relation of the SIPG, $\lambda_{ph,h}^s(\xi)$, satisfies: firstly, for $h$ fixed, $\lambda_{ph,h}^s(\xi)\to\lambda_h(\xi)$ as $s\to\infty$, for all $\xi\in\Pi_h$. On the other hand, for $s>1$ and $\xi\in\Pi_h$ fixed, $\lambda_{ph,h}^s(\xi)\to|\xi|$ as $h\to 0$, which is the dispersion relation of the continuous wave equation. Contrarily, as $s\to\infty$ for fixed $h$ or as $h\to 0$ for fixed $s>1$ and $\xi\in\Pi_h$, $\lambda_{sp,h}^s(\xi)\to+\infty$. Moreover, for all $s>1$, $\lambda_{ph,h}^s(\xi)$ is an increasing function of $\xi$ and $\omega_h(\xi)\leq\lambda_{ph,h}^s(\xi)\leq\lambda_h(\xi)$, for all $\xi\in\Pi_h$, and the physical and spurious group velocities, $\partial_{\xi}\lambda_{ph,h}^s(\xi)$ and $\partial_{\xi}\lambda_{sp,h}^s(\xi)$, have the following properties:
\begin{figure}[!]
\begin{center}
  \includegraphics[width=6cm, height=3.5cm]{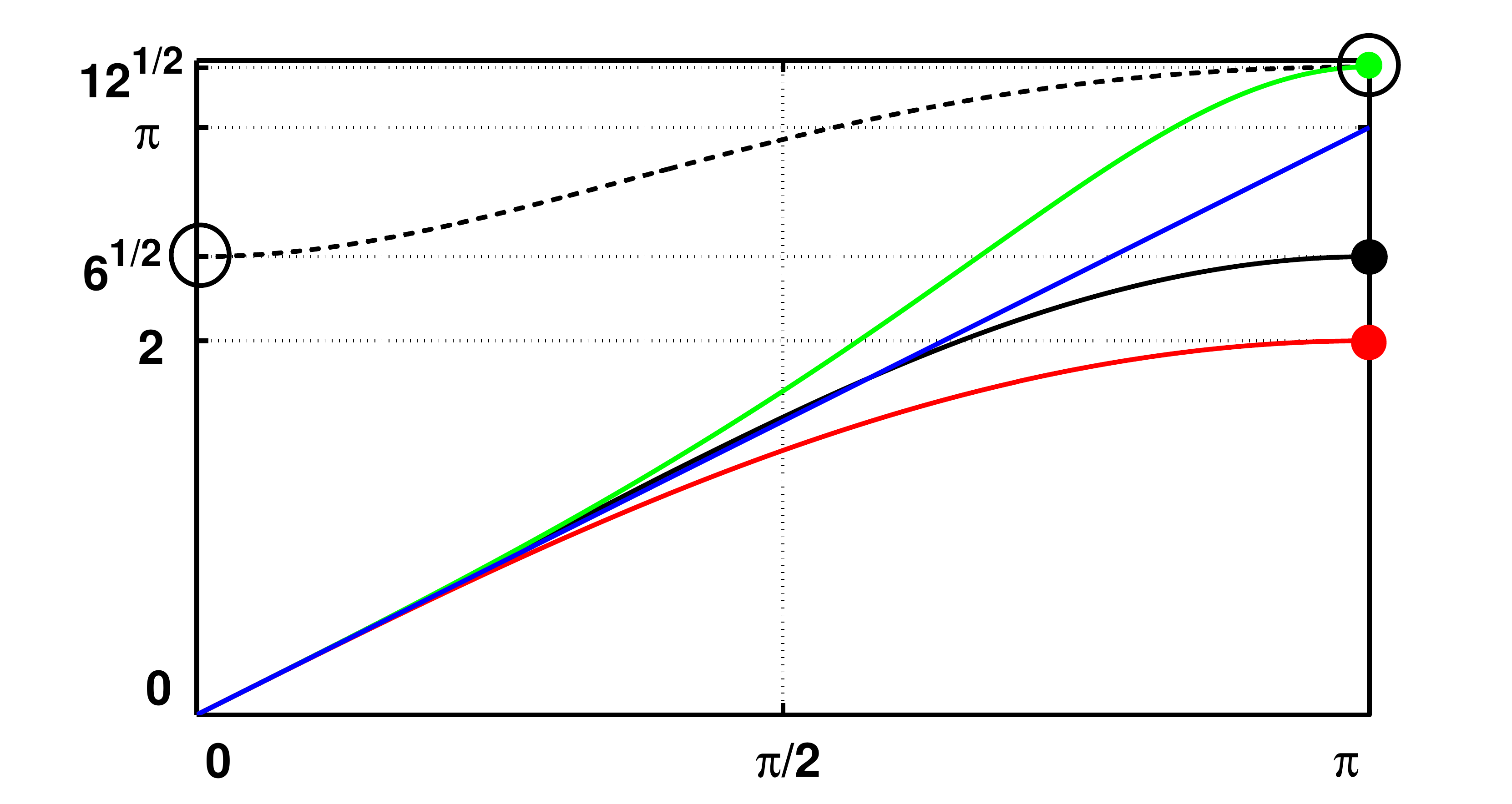}\includegraphics[width=6cm, height=3.5cm]{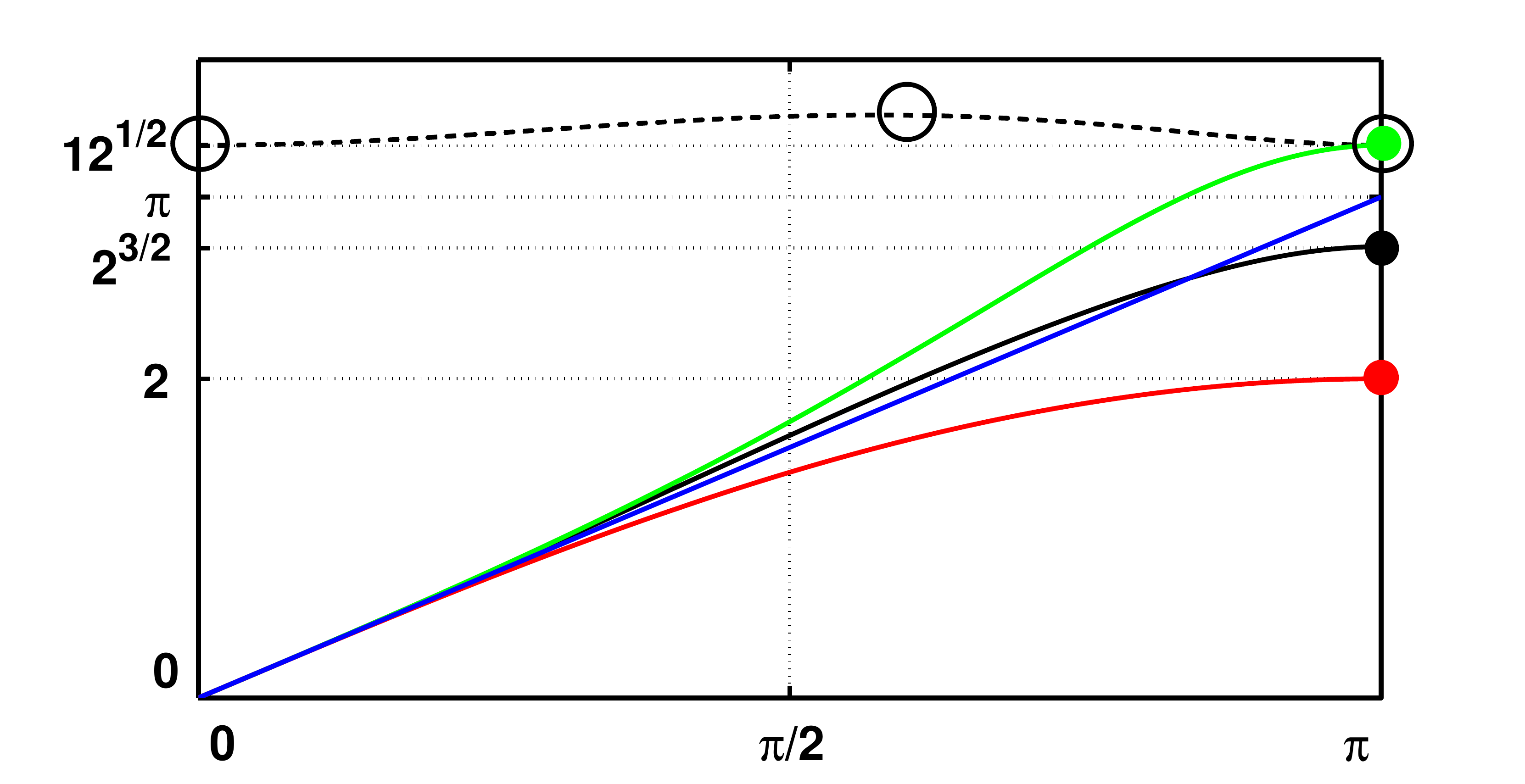}\\
  \includegraphics[width=6cm, height=3.5cm]{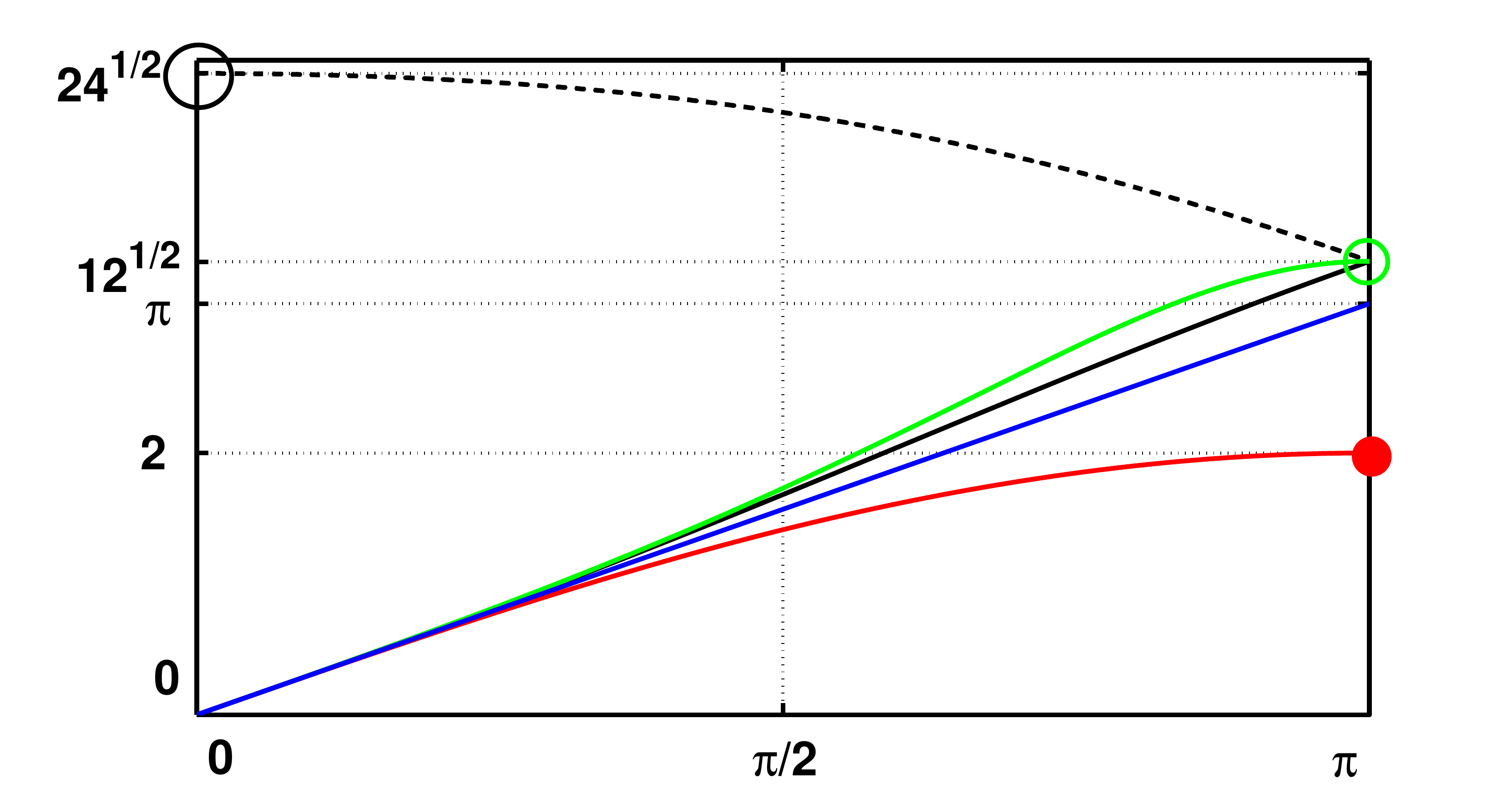}\includegraphics[width=6cm, height=3.5cm]{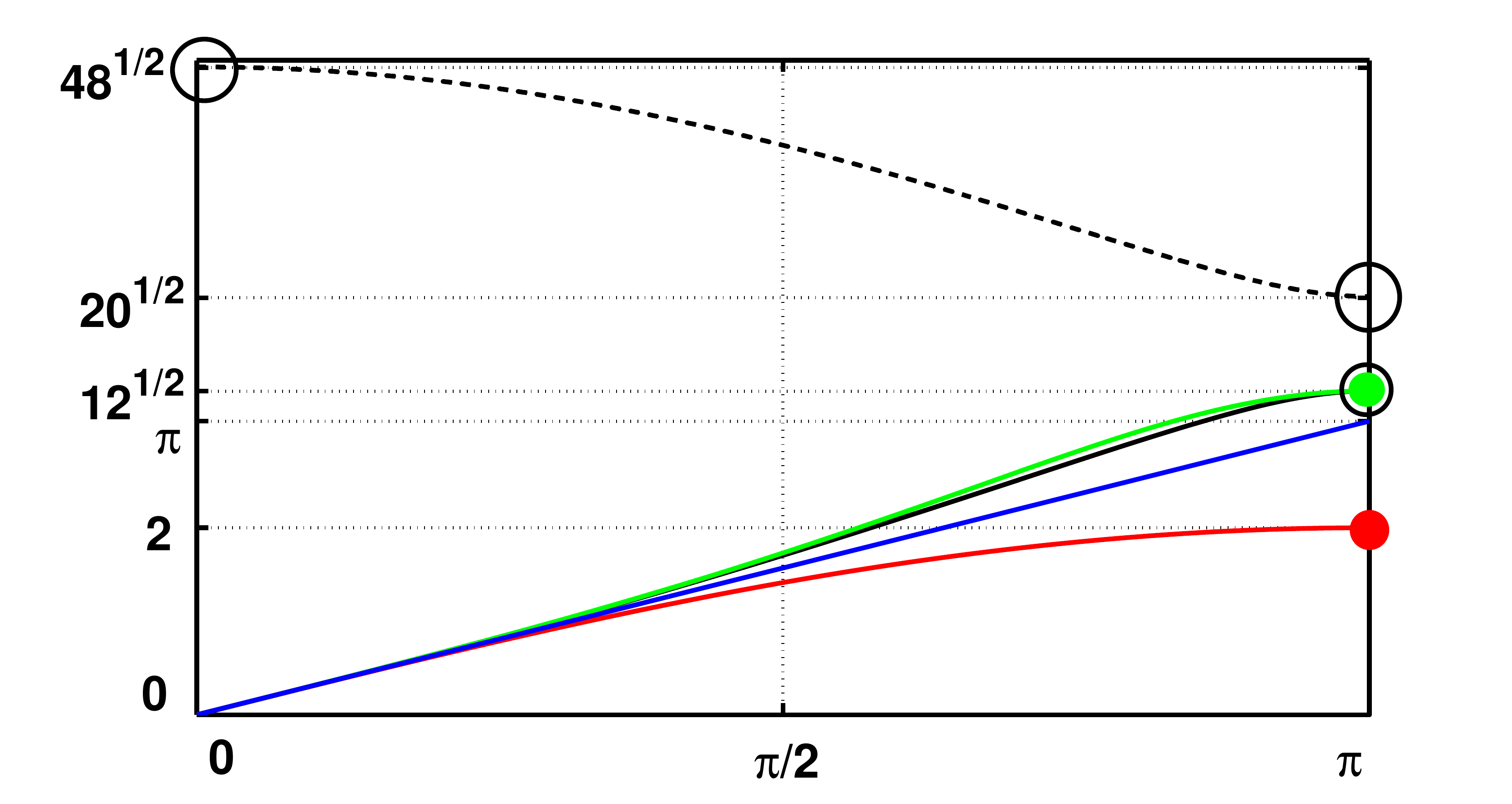}
\caption{Physical (black) and spurious (dotted black) dispersion relations, $\lambda_{ph,1}^s(\xi)$ and $\lambda_{sp,1}^s(\xi)$, for $s=1.5$ (top, left), $s=2$ (top, right), $s=3$ (bottom, left), $s=5$ (bottom, right) compared to the ones corresponding to the continuous wave equation (blue), $\xi$, and to its finite difference (red) and $P_1$-classical finite element (green) semi-discretizations, $\omega_1(\xi)$ and $\lambda_1(\xi)$. The marked points are wave numbers where the corresponding group velocities vanish.}\label{Dispersions}
  \end{center}
\end{figure}
\smallskip\begin{itemize}\item[1.] For all $s>1$, $\lim\limits_{\xi\to 0}\partial_{\xi}\lambda_{ph,h}^s(\xi)=1$ and $\lim\limits_{\xi\to 0}\partial_{\xi}\lambda_{sp,h}^s(\xi)=0$.
\item[2.] For all $s\in(1,\infty)\setminus\{3\}$, $\lim\limits_{\xi\to\pm\pi/h}\partial_{\xi}\lambda_{ph,h}^s(\xi)=\lim\limits_{\xi\to\pm\pi/h}\partial_{\xi}\lambda_{sp,h}^s(\xi)=0$.
\item[3.] $\lim\limits_{\xi\to\pm\pi/h}\partial_{\xi}\lambda_{ph,h}^3(\xi)=1$ and $\lim\limits_{\xi\to\pm\pi/h}\partial_{\xi}\lambda_{sp,h}^3(\xi)=-1$.
\item[4.] For all $s\in(1,\infty)$, the physical dispersion relation is increasing in $\xi$ for all $s>1$,  whereas the spurious one has several monotonicity ranges, according to the stabilization parameter $s$, as one can see in Figure \ref{Dispersions}.
\end{itemize}

\smallskip\noindent\textbf{2. Localized waves and observability inequalities.} For the continuous problem (\ref{ContWaveEqn}), it is well known that for any initial data $(u^0,u^1)\in \dot{H}^{1}(\rr)\times L^2(\rr)$ and any observability time $T>2$, there exists a constant $C(T)>0$ s.t. the following observability inequality holds (cf. \cite{ZuaUnbDom}):
\begin{equation}E(u^0,u^1)\leqslant C(T)\int\limits_0^TE_{\Omega}(u^0,u^1,t)\,dt,\label{ContObsIneq}\end{equation}
where $\Omega:=\rr\setminus(-1,1)$ and $E(u^0,u^1)$ and $E_{\Omega}(u^0,u^1,t)$ denote respectively the total energy (which is conserved in time) and the energy concentrated in $\Omega$ at time $t$, given explicitly by
\begin{equation}E(u^0,u^1)=\frac{1}{2}\int\limits_{\rr}(|\partial_tu(x,t)|^2+|\partial_xu(x,t)|^2)\,dx,\quad E_{\Omega}(u^0,u^1,t)=\frac{1}{2}\int\limits_{\Omega}(|\partial_tu(x,t)|^2+|\partial_xu(x,t)|^2)\,dx.\nonumber\end{equation}

The time $T^*=2$ is sharp, given by the so-called Geometric Control Condition
(GCC), requiring all rays of Geometric Optics to enter the observation region during the observability time. When the GCC does not hold, the observability property fails because of the Gaussian beam solutions localized around a
bi-characteristic ray escaping the observation region (see \cite{ErvZuaBook}).

We also analyze the SIPG version of the observability inequality:
\begin{equation}E_h^s(\overrightarrow{U}^{h,0},\overrightarrow{U}^{h,1})\leq C_h^s(T)\int\limits_0^TE_{\Omega,h}^s(\overrightarrow{U}^{h,0},\overrightarrow{U}^{h,1},t)\,dt.\label{DGObsIneq}\end{equation}
Here, $E_h^s(\overrightarrow{U}^{h,0},\overrightarrow{U}^{h,1})$ is the discrete total energy, conserved in time, given by
$$E_h^s(\overrightarrow{U}^{h,0},\overrightarrow{U}^{h,1})=\frac{1}{2}\Big(\langle R_h^s\overrightarrow{U}^{h,0},\overrightarrow{U}^{h,0}\rangle+
\langle M_h\overrightarrow{U}^{h,1},\overrightarrow{U}^{h,1}\rangle\Big),$$
where $\langle\cdot,\cdot\rangle$ is the inner product in $\ell^2$, and $E_{\Omega,h}^s(\overrightarrow{U}^{h,0},\overrightarrow{U}^{h,1},t)$ is the discrete energy concentrated in $\Omega$ at time $t$, defined as above but with $\langle\cdot,\cdot\rangle$ replaced by the local scalar product
$\langle\cdot,\cdot\rangle_{\Omega}$, which is the inner product in $\ell^2(\{j:x_j\in\Omega\})$. For all $T,h >0$, the inequality (\ref{DGObsIneq}) holds, with a finite constant $C^s_h(T)$. Our goal here is to analyze its behavior as $h \to 0$.

When the vector valued initial data $\widehat{U}^{h,i}$ in (\ref{DiscWaveFour}), $i=0,1$, are of the form
\begin{equation}\widehat{U}^{h,i}(\xi)=P_{ph,h}^s(\xi)\widehat{u}^{h,i}(\xi),\label{InitialDataPh}\end{equation}
the corresponding solutions of (\ref{DiscWaveFour}) involve only the physical dispersion relation:
\begin{equation}\widehat{U}^h(\xi,t)=P_{ph,h}^s(\xi)\frac{1}{2}\sum\limits_{\pm}\left(\widehat{u}^{h,0}(\xi)
\pm\frac{\widehat{u}^{h,1}(\xi)}{i\lambda_{ph,h}^s(\xi)}\right)\exp(\pm it\lambda_{ph,h}^s(\xi)).\label{SolnPh}\end{equation}
Considering solutions concentrated in wave packets and a stationary phase like argument
allow showing that, whatever $s$ and $T$ are, $C_h^s(T)$ blows up at an arbitrarily large polynomial rate as $h \to 0$:
\begin{proposition}Let $T>0$ be given with a semi-discrete ray $x_{ph}(t)=x^*-t\partial_{\xi}\lambda_{ph,1}^s(\eta_0)$ that does not enter the observation region in time $T$. Consider $\gamma:=\gamma(h)>0$ such that $\gamma>>1$ and $h\gamma<<1$. For $\phi\in\mathcal{S}(\rr)$, consider the semi-discrete wave equation (\ref{DiscWaveFour}) with initial data $\widehat{U}^{h,i}(\xi)$ satisfying (\ref{InitialDataPh}) with \begin{equation}\widehat{u}^{h,0}(\xi)=\sqrt{\frac{2\pi}{\gamma}}\widehat{\phi}\left(\frac{\xi-\xi_0}{\gamma}\right)
\exp(-ix^*(\xi-\xi_0)) \chi_{\Pi_h}(\xi)\mbox{ and }\ \widehat{u}^{h,1}(\xi)=i\lambda_{ph,h}^s(\xi)\widehat{u}^{h,0}(\xi).\label{initialdataSylvain}\end{equation}

Then for all $\alpha\in\rr_+$, the observability constant $C_h^s(T)$ in (\ref{DGObsIneq}) satisfies $C_h^s(T)\geqslant C_{\alpha}(\phi,T,s)\gamma^{\alpha}$.\end{proposition}
A similar result holds when the numerical solution is concentrated on the spurious mode.

\smallskip\noindent\textbf{3. Filtering mechanisms.} In what follows, we introduce a filtering mechanism aimed at recovering the uniformity as $h$ tends to zero of the observability constant $C_h^s(T)$ in (\ref{DGObsIneq})  within a subclass of solutions of the numerical approximation scheme. For $\delta\in(0,1)$, set $\Pi_h^{\delta}:=[-\pi\delta/h,\pi\delta/h]$ and let us define the space of Fourier filtered data $I_h^{\delta}=\{\overrightarrow{f}\in\ell^2(h\zz):\mbox{supp}(\widehat{f}^h)\subset \Pi_h^{\delta}\}$. It can be proved that if in (\ref{DiscWaveEqn}) we consider initial data $\overrightarrow{U}^{h,i}$, $i=0,1$, verifying the condition (\ref{InitialDataPh}) and such that $\widehat{u}^{h,i}\in I_h^{\delta}$ for $i=0,1$, then there exists an uniform time $T_{ph}^{s,\delta}$ such that, for all $T>T_{ph}^{s,\delta}$, the observability inequality (\ref{DGObsIneq}) holds uniformly as $h\to 0$. These data lead to solutions whose energy is concentrated on the low frequencies of the physical mode for which the group velocity of propagation is uniform.

However, our goal (as described in the pioneering work \cite{GloStokesWaves}) is to introduce a filtering mechanism that does not require the use of the Fourier transform, but that rather might be implemented in the numerical mesh, directly. The filtering mechanism we propose can be implemented in two steps. First, the initial data are taken so that their jumps vanish, i.e. $\overrightarrow{J}^{h,i}\equiv 0, i=0,1$.
Second, their average part, $\overrightarrow{A}^{h,i}$, is obtained by a bi-grid algorithm (analyzed in \cite{GloLioHeBook}), i.e. $A_{2j}^{i}=(A_{2j+1}^{i}+A_{2j-1}^{i})/2,\ \forall j\in\zz,\ \forall i=0,1$.

Although the corresponding solutions also excite the spurious spectral component, their energy is concentrated on the low frequency physical components. Thus, using the arguments in \cite{LivEZ}, one can show that, for $T$ large enough independent of $h$, the semi-discrete observability inequality (\ref{DGObsIneq}) holds uniformly in this class of filtered numerical solutions too. By duality, this implies a result of uniform  (with respect to $h$) controllability of a suitable projection of the solutions of the numerical approximation scheme.

\medskip\noindent\textbf{Acknowledgements.} Both authors were partially supported by the Grant MTM2008-03541 of the MICINN, Spain, and the ERC Advanced Grant FP7-246775 NUMERIWAVES.


\begin{thebibliography}{00}
\bibitem{AinsSisp} M.~Ainsworth, \newblock{\em Dispersive Behaviour of High Order Discontinuous Galerkin Finite Element Method}, Journal of Computational Physics, 198 (1) (2004), 106--130.


\bibitem{BuffaEig} P.~Antonietti, A.~Buffa, I.~Perugia, \newblock{\em Discontinuous Galerkin approximation of the Laplace eigenproblem}, Comput. Methods Appl. Mech. Engrg., 195 (25--28) (2006), 3483--3505.


\bibitem{BreUnif} D.N.~Arnold, F.~Brezzi, B.~Cockburn, L.D.~Marini, \newblock{\em Unified analysis of Discontinuous Galerkin Methods for Elliptic Problems}, SIAM J. Numer. Anal., 39 (2002), 1749--1779.


\bibitem{ErvZuaBook} S.~Ervedoza, E.~Zuazua, \newblock{\em Propagation, observation and numerical approximation of waves}, in preparation.

\bibitem{GloStokesWaves} R.~Glowinski, \newblock{\em Ensuring the well-posedness by analogy; Stokes problem and boundary control for the wave equation}, Journal of Computational Physics, 103 (2) (1992), 189--221.

\bibitem{GloLioHeBook} R.~Glowinski, J.-L.~Lions, J.~He, \newblock{\em Exact and approximate controllability for distributed parameter systems: a numerical approach}, Encyclopedia of Mathematics and its Applications, Cambridge University Press, 117 (2008).


\bibitem{LivEZ} L.~Ignat, E.~Zuazua, \newblock{\em Convergence of a multi-grid method for the control of waves}, J. Eur. Math. Soc., 11 (2009), 351–-391.

\bibitem{LivEZDispSchr} L.~Ignat, E.~Zuazua, \newblock{\em Numerical dispersive schemes for the nonlinear  Schr\"{o}dinger equation}, SIAM. J. Numer. Anal., 47(2) (2009), 1366--1390.

\bibitem{MarZuaCRASFD} A.~Marica, E.~Zuazua, \newblock{\em Localized solutions for the finite difference semi-discretization of the wave equation}, C.R. Acad. Sci. Paris, to appear.


\bibitem{ZuaUnbDom} E.~Zuazua, \newblock{\em Exponential decay for the semilinear wave equation with localized damping in unbounded domains}, J. Math. Pures Appl., 70 (1991), 513--529.

\bibitem{ZuaPOC} E.~Zuazua, \newblock{\em Propagation, Observation, Control and Numerical Approximations of Waves}, SIAM Review, 47(2)(2005), 197--243.

\end{thebibliography}
\end{document}